\documentclass[12pt]{article}

\usepackage{authblk}
\usepackage{color}
\usepackage[initials]{amsrefs}
\usepackage{amssymb}
\usepackage{amsmath}
\usepackage{amscd}
\usepackage{amsthm}
\usepackage[all]{xy}
\usepackage{nicefrac}
\usepackage{units}

\usepackage[mathscr]{eucal}
\usepackage{enumerate}

\def \Q {{\mathbb Q}}
\def \Z {{\mathbb Z}}

\def \R {{\mathbb R}}

\def \G {{\mathbb G}}

\newtheorem{theo}{Theorem}[section]
\newtheorem{thm}[theo]{Theorem}
\newtheorem{prop}[theo]{Proposition}
\newtheorem{lem}[theo]{Lemma}
\newtheorem{cor}[theo]{Corollary}
\newtheorem{defi}[theo]{Definition}
\newtheorem{exa}[theo]{Example}

\def\beq{\begin{equation} \label}

\newtheorem*{claim*}{Claim}
\newtheorem*{exa*}{Example}
\newtheorem*{rem*}{Remark}
\newtheorem*{rems*}{Remarks}
\newtheorem*{fact*}{Fact}
\newtheorem*{Def*}{Definition}

\newcommand{\bthe}{\begin{thm}}
\newcommand{\ble}{\begin{lem}}
\newcommand{\bpr}{\begin{prop}}
\newcommand{\bco}{\begin{cor}}
\newcommand{\bde}{\begin{defi}}
\newcommand{\ethe}{\end{theo}}
\newcommand{\ele}{\end{lem}}
\newcommand{\epr}{\end{prop}}
\newcommand{\eco}{\end{cor}}
\newcommand{\ede}{\end{defi}}

\def \Gal {{\rm{Gal}}}

\def \Ker {{\rm Ker}}

\DeclareFontFamily{U}{wncy}{}
\DeclareFontShape{U}{wncy}{m}{n}{%
<5>wncyr5%
<6>wncyr6%
<7>wncyr7%
<8>wncyr8%
<9>wncyr9%
<10>wncyr10%
<11>wncyr10%
<12>wncyr6%
<14>wncyr7%
<17>wncyr8%
<20>wncyr10%
<25>wncyr10}{} \DeclareMathAlphabet{\cyr}{U}{wncy}{m}{n}
\begin{document}

\title{On the equation $x^2-Dy^2=n$}
\author{Dasheng Wei}

\affil{Academy of Mathematics and System Science, CAS, Beijing
100190, P.R.China\\email:dshwei@amss.ac.cn}

\date{\today}

\maketitle

 \bigskip
\section*{\it Abstract}

We propose a method to determine the solvability of the diophantine
equation $x^2-Dy^2=n$ for the following two cases:

$(1)$ $D=pq$, where $p,q\equiv 1 \mod 4$ are distinct primes with
$\left(\frac{q}{p}\right)=1$ and
$\left(\frac{p}{q}\right)_4\left(\frac{q}{p}\right)_4=-1$.

$(2)$ $D=2p_1p_2\cdots p_m$, where $p_i\equiv 1 \mod 8,1\leq i\leq
m$ are distinct primes and $D=r^2+s^2$ with $r,s \equiv \pm 3 \mod
8$.

\bigskip

{\it MSC classification} : 11D09; 11E12


\bigskip

{\it Keywords} : quadratic form, torus, Hilbert class field,
reciprocity law.

\section*{Introduction} \label{sec.notation}
Let $D$ be a non-square integer. The question of whether the
equation
\begin{equation}\label{bi} x^2-Dy^2=n, n\in \Z\end{equation}
has an integral solution is a very old one (see \cite{Di}). We may
recast the question in the language of algebraic geometry and ask
whether the affine scheme over $\Z$ defined by (\ref{bi}) has an
integral point. It's well-known that the generic fiber of this
affine scheme is a principal homogenous space of tori when $n\neq
0$. Recently, Harari \cite{Ha08} showed that the Brauer-Manin
obstruction is the only obstruction for existence of the integral
points of such scheme. Fei Xu and the author gave another proof of
the result in \cite{WX} and \cite{WX2}. In this paper we consider
the solvability of (\ref{bi}) by using the method in \cite{WX}.

It should be pointed out that the method in \cite{WX} only produces
the idelic class groups of $\Q(\sqrt{D})$ and these idelic class
groups are not unique. In order to get the explicit conditions for
the solvability, one needs further to construct the explicit abelian
extensions of $\Q(\sqrt{D})$ corresponding to the idelic class
groups. Such explicit construction is a wide open problem in
general.

Notation and terminology are standard if not explained. Let
$E=\Q(\sqrt{D})$ and $\frak o_E$ be the ring of integers of $E$,
$\Omega_E$ the set of places of $E$ and $\infty$ the set of infinite
places of $E$. Let $E_\frak p$ be the completion of $E$ at $\frak p$
and $\frak o_{E_\frak p}$ the local completion of $\frak o_E$ at
$\frak p$ for each $\frak p\in \Omega_E$. Write $\frak o_{E_\frak
p}=E_\frak p$ for $\frak p\in \infty$ and denote the adele ring
(resp. the idele group) of $E$ by $\Bbb A_E$ (resp. $\Bbb I_E$).

Let $\bold X_{(n)}$ denote the affine scheme over $\frak o_{F}$
defined by $x^2-Dy^2=n$ for a non-zero integer $n$. Let
$X_{(n)}=\bold X_{(n)} \times_{\Z} {\Q}$. Obviously $f=x+y\sqrt{D}$
is an invertible function on $X_{(n)}\otimes_\Q E$. And $f$ induces
a natural map
$$f_E:\ \ \ X_{(n)}(\Bbb A_\Q)\rightarrow \Bbb I_E.$$ The
restriction to $X_{(n)}(\Q_p)$ of $f_E$ can be written by
$$f_E[(x_p,y_p)]= \begin{cases} (x_p+y_p \sqrt{D}, x_p-y_p\sqrt{D}) \ \ \
& \text{if $p$ splits in $E/\Q$,} \\
x_p+y_p\sqrt{D} \ \ \ & \text{otherwise.}
\end{cases} $$
\begin{Def*}
Let $K_1,\cdots,K_m$ be finite abelian extensions over $E$. Let
$$\psi_{K_i/E}: \Bbb I_E\rightarrow \Gal(K_i/E) \text{ for } 1\leq
i\leq m$$ be the Artin map. We say that $n$ satisfies the Artin
condition of $K_1,\cdots,K_m$ if there is
$$\prod_{p\leq \infty}(x_p,y_p)\in \prod_{p\leq \infty}
\bold X_{(n)}(\Z_p)$$ such that
$$\psi_{K_i/E}(f_E[\prod_{p\leq \infty}(x_p,y_p)])=1_{i} \text{ for } i=1,\cdots,m$$
where $1_i$ is the identity element of $\Gal(K_i/E)$.
\end{Def*}

By the class field theory, it is a necessary condition for $\bold
X_{(n)}(\Z) \neq \emptyset$ that $n$ satisfies the Artin condition
of $K_1,\cdots,K_m$. And there is a finite abelian extension $K/E$
that is independent on $n$, such that the Artin condition of $K$ is
also sufficient for $\bold X_{(n)}(\Z) \neq \emptyset$ (Corollary
2.8 in \cite{WX}). For example, let $L=\Z+\Z\sqrt{D}$ and $H_L$ be
the ring class field corresponding to the order $L$, then the Artin
condition of $H_L$ is sufficient for $\bold X_{(n)}(\Z) \neq
\emptyset$ if $D<0$ (Proposition 3.1 in \cite{WX}). However, the
Artin condition of $H_L$ is not always sufficient for general $D$.

Let $E=\Q(\sqrt{D})$ and let $T$ be the torus $$R^1_{E/F}(\Bbb
G_m)=\Ker[R_{E/F}(\Bbb G_{m,E})\rightarrow \Bbb G_{m,\Q}],$$ here
$R_{E/F}$ denotes the Weil's restriction (see \cite{Milne98}, p.
225). Denote $\lambda$ to be the embedding from $T$ to $
R_{E/F}(\Bbb G_{m,E})$. Obviously $\lambda$ induces a natural group
homomorphism $\lambda_E:T(\Bbb A_F)\rightarrow \Bbb I_E.$ Let $\bold
T$ be the group scheme over $\Z$ defined by $x^2-Dy^2=1$ and
$T=\bold T \times_{\Z} \Q$. The generic fiber of $\bold X_{(n)}$ is
a principal homogenous space of the torus $T$. Since $\bold T$ is
separated over $\Z$, we can view $\bold T(\Z_p)$ as a subgroup of
$T(\Q_p)$. The following result can be founded in \cite{WX}
(Corollary 2.20).
\begin{prop} \label{multiple} Let $K_1/E$ and $K_2/E$ be finite abelian
extensions
 such that the group homomorphism induced by
$\lambda_E$ (we also denote it by $\lambda_E$)
$$\lambda_E: \nicefrac{T(\Bbb A_\Q)}{T(\Q)\prod_{p\leq \infty}\bold T(\Z_p)}
\longrightarrow \nicefrac{\Bbb I_{E}}{E^*N_{K_1/E}(\Bbb
I_{K_1})}\times \nicefrac{\Bbb I_{E}}{E^*N_{K_2/E}(\Bbb I_{K_2})}
$$ is well-defined and injective, where well-defined means
$$ \lambda_E(T(\Q)\prod_{p\leq \infty}\bold T(\Z_p))\subset (E^*N_{K_1/E}(\Bbb I_{K_1}))\cap (E^*N_{K_2/E}(\Bbb I_{K_2})).$$ Then $\bold X_{(n)}(\Z)\neq \emptyset$ if and only if
$n$ satisfies the Artin condition of $K_1$ and $K_2$.

\end{prop}


Let $p$ and $q$ are distinct primes. The following facts are
well-know:

 (1) If $p\equiv 3 \mod
4$ or $q\equiv 3 \mod 4$, then $x^2-pqy^2=-1$ is not solvable over
$\Z_2$.

 (2) If $p$ and $q$ are of the form $4k+1$ and
$\left(\frac{p}{q}\right)=-1$, then $x^2-pqy^2=-1$ is solvable over
$\Z$ (\cite{D}, p. 228).

 (3) If $p$ and $q$ are of the form $4k+1$
with $\left(\frac{p}{q}\right)=1$ and
$\left(\frac{p}{q}\right)_4=\left(\frac{q}{p}\right)_4=-1$, then
$x^2-pqy^2=-1$ is also solvable over $\Z$ (\cite{D}, p. 228).

For the above three cases, the equation
\begin{equation}\label{equ2}x^2-pqy^2=n\end{equation} is solvable over $\Z$ if and only if
$n$ satisfies the Artin condition of $H_L$ by Proposition 4.1 of
\cite{WX}. Therefore we only need to consider the solvability of
(\ref{equ2}) when $p,q \equiv 1 \mod 4 $,
$\left(\frac{p}{q}\right)=1$ and $\left(\frac{p}{q}\right)_4=1\text{
or }\left(\frac{q}{p}\right)_4=1$. In \S 1, Theorem \ref{0.1}, we
consider the solvability of (\ref{equ2}) when $p,q \equiv 1 \mod 4
$, $\left(\frac{q}{p}\right)=1$ and
$\left(\frac{p}{q}\right)_4\left(\frac{q}{p}\right)_4=-1$. As an
application, we reprove Scholz and Brown's result  (\cite{Sc},
\cite{Br}) about solvability  of the equations $x^2-pqy^2=-1,p,q$
(see Corollary \ref{Sc-Br}).



Let $d=p_1p_2\cdots p_m$ where $p_i\equiv 1 \mod 8, 1\leq i \leq m$
are distinct primes. In \S 2, Theorem \ref{0.2}, we consider the
solvability of $x^2-2dy^2=n$. As application of Theorem \ref{0.2},
the solvability of the equations $x^2-2dy^2=-1,\pm 2$ are considered
in \S 3. In particular, we recover the main result of Pall in
\cite{Pa} (see the remarks following Proposition \ref{9-case} and
Corollary \ref{pall} below).


\section{The solvability of $x^2-pqy^2=n$}
Let $p$ and $q$ be distinct primes of the form $4k+1$ with
$\left(\frac{q}{p}\right)=1$.  Then the equation $x^2-pqy^2=pz^2$ is
solvable over $\Z$ by the Hasse principle. Fix an integral solution
$(x_0,y_0,z_0)$ of the equation such that $x_0>0$ and $(x_0,y_0)=1$.
Let $E=\Q(\sqrt{pq})$ and $\Theta=E(\sqrt{x_0-y_0\sqrt{pq}})$. Then
$\Theta$ is totally real and $\Theta/E$ is unramified over all
primes except the primes over $2p$.

\begin{lem} \label{computation-pq} Let $p$ and $q$ be distinct primes of the form $4k+1$  with
$\left(\frac{q}{p}\right)=1$. Let $l=2 \text{ or } p$. If $x_l$ and
$y_l$ in $\Z_l$ satisfy $x_l^2-pqy_l^2=1$, then the quadratic
Hilbert symbol
$$\prod_{v\mid l}\left(\frac{x_l-y_l\sqrt{pq}, x_0-y_0\sqrt{pq}}{v}\right)=1$$
where $v\in \Omega_E$.
\end{lem}
\begin{proof} $(1)$ Assume $l$ is split in  $E/\Q$. Then $l=2$ since $E/\Q$ is ramified at $p$.
Then we have $$\aligned &\prod_{v\mid
2}\left(\frac{x_2-y_2\sqrt{pq},
x_0-y_0\sqrt{pq}}{v}\right)\\
&=\left(\frac{x_2-y_2\sqrt{pq},
x_0-y_0\sqrt{pq}}{2}\right)\cdot\left(\frac{x_2+y_2\sqrt{pq},
x_0+y_0\sqrt{pq}}{2}\right)\\
&=\left(\frac{x_2-y_2\sqrt{pq},
x_0-y_0\sqrt{pq}}{2}\right)\cdot\left(\frac{(x_2-y_2\sqrt{pq})^{-1},
x_0+y_0\sqrt{pq}}{2}\right)\\
&=\left(\frac{x_2-y_2\sqrt{pq},
x_0-y_0\sqrt{pq}}{2}\right)\cdot\left(\frac{x_2-y_2\sqrt{pq},
x_0+y_0\sqrt{pq}}{2}\right)\\
&=\left(\frac{x_2-y_2\sqrt{pq}, p}{2}\right)=1,
\endaligned
$$
the last equation holds since $x_l-y_l\sqrt{pq}$ is a unit in $\Z_2$
and $p\equiv 1 \mod 4$.

$(2)$ Assume $l$ is not split in $E/\Q$. Let $\xi=x_l+y_l\sqrt{pq}$.
Let $v$ be the unique place of $E$ over $l$. Since $N_{E_v/\Bbb
Q_l}(\xi)=1$, there exists $\alpha\in E_v^*$ such that
$\xi=\sigma(\alpha)\alpha^{-1}$ by Hilbert's theorem 90, where
$\sigma$ is the non-trivial element in $Gal(E_v/\Bbb Q_l)$. Then
$$\aligned &\left(\frac{\xi, x_0-y_0\sqrt{pq}}{v}\right)=\left(\frac{\sigma(\alpha)\alpha^{-1},
x_0-y_0\sqrt{pq}}{v}\right)\\
=&\left(\frac{N_{E_v/\Bbb Q_l}(\alpha),
x_0-y_0\sqrt{pq}}{v}\right)=\left(\frac{N_{E_v/\Bbb
Q_l}(\alpha),p}{l}\right).
\endaligned
$$
Suppose $l=2$. Since $E_v/\Q_2$ is unramified and $p\equiv 1 \mod
4$, one has $$\left(\frac{\xi, x_0-y_0\sqrt{pq}}{v}\right) =
\left(\frac{N_{E_v/\Bbb Q_2}(\alpha),p}{2}\right)=1.$$ Suppose
$l=p$. Since $\left(\frac{q}{p}\right)=1$ and $E_v=\Bbb
Q_p(\sqrt{pq})$, one has
$$\left(\frac{\xi, x_0-y_0\sqrt{pq}}{v}\right) =
\left(\frac{N_{E_v/\Bbb
Q_p}(\alpha),p}{p}\right)=\left(\frac{N_{E_v/\Bbb
Q_p}(\alpha),pq}{p}\right)=1.$$
\end{proof}

We make use of the following interesting result, due to K. Burde
(see \cite{B}, p. 183):
\begin{lem} \label{Burde} If $p=a^2+b^2,q=c^2+d^2,a\equiv c\equiv 1,b\equiv d\equiv 0 (\text{mod }2), a,b,c,d>0$, p and q are distinct primes, and $\left(\frac{q}{p}\right)=1$, then
$$\left(\frac{p}{q}\right)_4\left(\frac{q}{p}\right)_4=(-1)^{(p-1)/4}\left(\frac{ad-bc}{p}\right).$$
\end{lem}

\begin{lem} \label{computation2-pq}
Let $p$ and $q$ be distinct primes of the form $4k+1$. Suppose
$\left(\frac{q}{p}\right)=1$ and
$\left(\frac{p}{q}\right)_4\left(\frac{q}{p}\right)_4=-1$. Let $l=2
\text{ or } p$ and let $x_l$ and $y_l$ in $\Z_l$ satisfy
$x_l^2-pqy_l^2=-1$. Then the quadratic Hilbert symbol
$$\prod_{l\mid 2p}\prod_{v\mid l}\left(\frac{x_l+y_l\sqrt{pq}, x_0-y_0\sqrt{pq}}{v}\right)=-1$$
where $v\in \Omega_E$.
\end{lem}
\begin{proof} Let $p=a^2+b^2,q=c^2+d^2,a\equiv c\equiv 1,b\equiv d\equiv 0 \text{ (mod } 2), a,b,c,d>0$. And let $r=ad-bc, s=ac+bd$. Then $r^2+s^2=pq$ and $s$ is odd.
Obviously $(r/s)^2\equiv -1 \mod p$. So
$$\begin{aligned}\left(\frac{s}{p}\right)&=\left(\frac{r}{p}\right)\left(\frac{-1}{p}\right)_4
=\left(\frac{r}{p}\right)\left(\frac{2}{p}\right)
=\left(\frac{r}{p}\right)(-1)^{(p-1)/4}\\
&=\left(\frac{p}{q}\right)_4\left(\frac{q}{p}\right)_4=-1
\end{aligned}$$ by Lemma \ref{Burde}.

By Lemma \ref{computation-pq}, we have
$$\aligned & \prod_{l\mid 2p}\prod_{v\mid l}\left(\frac{x_l+y_l\sqrt{pq}, x_0-y_0\sqrt{pq}}{v}\right)=\prod_{v\mid 2p}\left(\frac{(r-\sqrt{pq})s^{-1}, x_0-y_0\sqrt{pq}}{v}\right) \\
      &  = \prod_{v\mid 2p}\left(\frac{r-\sqrt{pq}, x_0-y_0\sqrt{pq}}{v}\right)\cdot \prod_{v\mid 2p}\left(\frac{s, x_0-y_0\sqrt{pq}}{v}\right)  \\
      & = \prod_{v\mid 2p}\left(\frac{r-\sqrt{pq}, x_0-y_0\sqrt{pq}}{v}\right)\cdot\left(\frac{s,pz_0^2}{2} \right)\cdot \left(\frac{s,pz_0^2}{p}
       \right)\\
      & = -\prod_{v\mid 2p}\left(\frac{r-\sqrt{pq}, x_0-y_0\sqrt{pq}}{v}\right). \endaligned  $$
By the Hilbert reciprocity law, one has
$$\prod_{v\mid 2p}\left(\frac{r-\sqrt{pq}, x_0-y_0\sqrt{pq}}{v}\right)
=\prod_{\frak p\nmid 2p}\left(\frac{r-\sqrt{pq},
x_0-y_0\sqrt{pq}}{\frak p}\right).$$ Since
$$(r-\sqrt{pq})(r+\sqrt{pq})=-s^2 \ \ \ \text{and} \ \ \
(x_0-y_0\sqrt{pq})(x_0+y_0\sqrt{pq})=pz_0^2$$ with $(r,s)=1$ and
$(x_0,y_0)=1$ respectively, one has $$ord_\frak p(r-\sqrt{pq})\equiv
ord_\frak p(x_0-y_0\sqrt{pq})\equiv 0 \mod 2
$$ for $\frak p\nmid 2p$ and $\frak p<\infty_E$. Since
$x_0-y_0\sqrt{pq}>0$ over $\frak p\in \infty_E$, one obtains that
$$\left(\frac{r-\sqrt{pq}, x_0-y_0\sqrt{pq}}{\frak p}\right)=1 $$ for all $\frak
p\nmid 2p$. One concludes that $$\prod_{l\mid 2p}\prod_{v\mid
l}\left(\frac{x_l+y_l\sqrt{pq}, x_0-y_0\sqrt{pq}}{v}\right)=-1.$$
\end{proof}

The following result was first proved by Scholz (see \cite{Sc}) and
was reproved by Brown (see \cite{Br}). Use our method, now we can
give a new proof.
\begin{cor}\label{Sc-Br} Let $p$ and $q$ be distinct primes of the form $4k+1$. Suppose
$\left(\frac{q}{p}\right)=1$ and
$\left(\frac{p}{q}\right)_4\left(\frac{q}{p}\right)_4=-1$. Then the
equation $x^2-pqy^2=-1$ is not solvable over $\Z$; the equation
$x^2-pqy^2=p$ is solvable over $\Z$ if and only if
$\left(\frac{q}{p}\right)_4=1$.
\end{cor}
\begin{proof}By Lemma \ref{computation2-pq}, one has $$\prod_{l\mid 2p}\prod_{v\mid
l}\left(\frac{x_l+y_l\sqrt{pq}, x_0-y_0\sqrt{pq}}{v}\right)=-1$$ for
$(x_l,y_l)\in \Z_l\times \Z_l$ with $x_l^2-pqy_l^2=-1$. This implies
that
$$ \psi_{\Theta/E}( f_E[\prod_{p\leq
\infty}(x_p,y_p)])=-1
$$ for any $\prod_{l\leq \infty}(x_l,y_l)\in \prod_{l\leq \infty}
\bold X_{(-1)}(\Bbb Z_l)$, since $\Theta/E$ is unramified over each
prime $v$ except $v\mid 2p$. Then $x^2-pqy^2=-1$ is not solvable
over $\Z$ by the class field theory.

One and only one of the three equations
$$x^2-pqy^2=-1,x^2-pqy^2=p,x^2-pqy^2=q$$ is solvable over $\Z$ (\cite{D},
p. 228). If $x^2-pqy^2=p$ is solvable over $\Z$, we know
$\left(\frac{q}{p}\right)_4=1$ (\cite{D}, p. 230). If
$\left(\frac{q}{p}\right)_4=1$, then
$\left(\frac{p}{q}\right)_4=-1$. So $x^2-pqy^2=q$ is not solvable
over $\Z$. We already know $x^2-pqy^2=-1$ is not solvable over $\Z$.
Therefore $x^2-pqy^2=p$ is solvable over $\Z$.
\end{proof}

By Lemma \ref{computation-pq} and \ref{computation2-pq}, now we can
give the main result of this section.
\begin{thm}\label{0.1}
Let $p$ and $q$ be distinct primes of the form $4k+1$. Suppose
$\left(\frac{q}{p}\right)=1$ and
$\left(\frac{p}{q}\right)_4\left(\frac{q}{p}\right)_4=-1$.  Then the
diophantine equation $x^2-pqy^2=n$ is solvable over $\Z$ if and only
if $n$ satisfies the Artin condition of $H_L$ and $\Theta$, where
$H_L$ is the ring class field corresponding to $L=\Z[\sqrt{pq}]$ and
$\Theta/E$ is a quadratic extension defined as above.
\end{thm}
\begin{proof}
Let $L=\Bbb Z + \Bbb Z \sqrt{pq}$. For any prime $l$,
$E_l=E\otimes_\Bbb Q \Bbb Q_l$ and $L_l$ is the $l$-adic completion
of $L$ inside $E_l$. Recall $T=R^1_{E/F}(\G_{m,E})$ and $\bold T$ is
the group scheme over $\Z$ defined by $x^2-pqy^2=1$, we have
$$T(\Q)=\{\xi \in E^*: \ N_{E/\Q}(\xi)=1\}$$ and
$$\bold T(\Z_p)=\{\xi \in L_ p^\times: \ N_{E_p/\Q_p}(\xi)=1\}.$$
And $L_\infty ^\times=E_\infty^*= \R^*\times \R^*$.

Let $l=2 \text{ or } p$. By Lemma \ref{computation-pq}, one has
$$\prod_{v\mid l}\left(\frac{\xi, x_0-y_0\sqrt{pq}}{v}\right)=1$$
for $(x_l,y_l)\in \Z_l\times \Z_l$ with $x_l^2-pqy_l^2=1$, where
$v\in \Omega_E$. This implies that
$$\lambda_E(\bold T(\Z_l))\subseteq E^* N_{\Theta/E}(\Bbb I_\Theta). $$
Since $\Theta/E$ is unramified over each prime $v$ except $v\mid
2p$, the natural group homomorphism
$$\lambda_E: \nicefrac{T(\Bbb A_\Q)}{T(\Q)\prod_{l\leq \infty}\bold T(\Z_l)}
\longrightarrow [\nicefrac{\Bbb I_E}{E^*  N_{\Theta/E}(\Bbb
I_{\Theta})}] \times [\nicefrac{\Bbb I_E}{E^* \prod_{l\leq \infty}
L_l^\times}]$$ is well-defined. By Proposition \ref{multiple}, we
only need to show $\lambda_E$ is injective.

Let $u\in \ker \lambda_E$. Then there are $\alpha\in E^*$ and $i\in
\prod_{l\leq \infty} L_l^\times$ with $\lambda_E(u)=\alpha i$. We
have
$$N_{E/\Q}(\alpha)=N_{E/\Q}(i)^{-1} \in \Q^* \cap
(\prod_{l\leq \infty}\Z_l^\times )=\{\pm 1\}. $$ If
$N_{E/\Q}(\alpha)\neq 1$, one obtains
$N_{E/\Q}(\alpha)=N_{E/\Q}(i)=-1$. Write $i=(i_v)_{v}\in \Bbb I_E$.
Since $\Theta/E$ is unramified over each prime $v$ except $v\mid
2p$, one concludes that $\psi_{\Theta/E}(i_v)$ is trivial for all
primes $v\nmid 2p$, where $i_v$ is regarded as an idele whose
$v$-component is $i_v$ and 1 otherwise. One gets
$$\psi_{\Theta/E}(\alpha i)=\psi_{\Theta/E}(i)=\prod_{v\mid 2p}\psi_{\Theta/E}(i_{v})=-1$$
by Lemma \ref{computation2-pq}, where $\psi_{\Theta/E}:\Bbb
I_E\rightarrow Gal(\Theta/E)$ is the Artin map. This contradicts to
$u \in \ker \lambda_E$.

Therefore $N_{E/\Q}(\alpha)=1$, one concludes that
$$ N_{E/\Q}(\alpha) =N_{E/\Q}(i) =1 \ \ \ \Rightarrow \ \ \ \alpha\in
T(\Q) \ \ \ \text{and} \ \ \ i\in \prod_{l\leq \infty}\bold
T(\Z_l).$$ So $\alpha i\in T(\Q)\prod_{l\leq \infty}\bold T(\Z_l)$.
\end{proof}

Finally we will use Theorem \ref{0.1} to give an explicit example.
For any integer $n$, one can write $n=(-1)^{s_0}
2^{s_1}13^{s_2}17^{s_3}{p_1}^{e_1}\cdots {p_g}^{e_g}$ and
$P(n)=\{p_1, \cdots, p_g \}$. Denote
$$\aligned & P_1=\{p\in P(n) : \ \left(\frac{13}{p}\right)=\left(\frac{17}{p}\right)=-1  \} \text{ and } P_2=\{p\in P(n): \
\left(\frac{221}{p}\right)=-1 \} \cr
 & P_3= \{p\in P(n): \
\left(\frac{13}{p}\right)=\left(\frac{17}{p}\right)=1 \ \text{and} \
x^4-238x^2+17 \equiv 0 \mod \ p \ \ \text{is solvable} \}.
\endaligned $$ Let $$n_1=\prod_{p_i\in
D(n)\setminus D_2} p_i^{e_i}$$

\begin{exa} \label{pell} Let $n$ be an integer with the above notation.
Then the equation $$x^2-221y^2=n$$ is solvable over $\Bbb Z$ if and
only if

(1)  $s_1$ is even, $(\frac{n_1}{17})=1$ and $(\frac{221}{p_i})=1$
for odd $e_i$.

(2)  $P_1\neq \emptyset$; or $$\prod_{p_i\in P_3}
(-1)^{e_i}\prod_{p_i\in P(n)\setminus P_2}
\left(\frac{-1}{p_i}\right)^{e_i}=(-1)^{s_0+s_2}\cdot
\left(\frac{n_1}{17}\right)_4$$ for $P_1=\emptyset$.
\end{exa}
\begin{proof} Since $\left(\frac{17}{13}\right)_4=-1$ and $\left(\frac{13}{17}\right)_4=1$, Theorem
\ref{0.1} can be applied. The ring class field $H_L$ associated to
the subring $\Z[\sqrt{221}]$ is $\Bbb Q(\sqrt{13},\sqrt{17})$. Since
the equation $x^2-221y^2=17$ has an integral solution for $x=119$
and $y=8$, one can choose $\Theta=E(\sqrt{119-8\sqrt{221}})$. We can
get this result by some computations.

\end{proof}

\section{The solvability of $x^2-2dy^2=n$}

Let $d=p_1\cdots p_m$ where $p_i\equiv 1 \mod 8, 1\leq i \leq m$ are
distinct primes. The equation $x^2-2dy^2=2z^2$ is solvable over $\Z$
by the Hasse principle. Fix an integral solution $(x_0,y_0,z_0)$ of
the equation such that $x_0>0$ and $(x_0,y_0)=1$. Let
$\Theta=E(\sqrt{x_0-y_0\sqrt{2d}})$. Then $\Theta$ is totally real
and $\Theta/E$ is unramified over all primes except the prime above
$2$ and $2$ is totally ramified in $\Theta/\Bbb Q$. First the
following lemmas will be proved.

\begin{lem} \label{computation} Let $d=p_1p_2\cdots p_m$
where $p_i\equiv 1 \mod 8, 1\leq i \leq m$ are distinct primes. Then
the quadratic Hilbert symbol
$$\left(\frac{\xi, x_0-y_0\sqrt{2d}}{v}\right)=1$$
for any $\xi\in E_v^* $ with $N_{E_v/\Bbb Q_2}(\xi)=1$, where $v$ is
the unique place of $E$ over $2$ and $(x_0, y_0)$ is given as above.
\end{lem}
\begin{proof}
Since $N_{E_v/\Bbb Q_2}(\xi)=1$, there exists $\alpha\in E_v^*$ such
that $\xi=\sigma(\alpha)\alpha^{-1}$ by Hilbert's theorem 90, where
$\sigma$ is the non-trivial element in $\Gal(E_v/\Bbb Q_2)$. Then
$$\aligned & \left(\frac{\xi, x_0-y_0\sqrt{2d}}{v}\right)=\left(\frac{\sigma(\alpha)\alpha^{-1},
x_0-y_0\sqrt{2d}}{v}\right)\\
= \ & \left(\frac{N_{E_v/\Bbb Q_2}(\alpha),
x_0-y_0\sqrt{2d}}{v}\right)=\left(\frac{N_{E_v/\Bbb
Q_2}(\alpha),2}{2}\right).
\endaligned
$$ Since
$d\equiv 1 \mod 8$ and $E_v=\Bbb Q_2(\sqrt{2d})$, one has
$$\left(\frac{\xi, x_0-y_0\sqrt{2d}}{v}\right) =
\left(\frac{N_{E_v/\Bbb
Q_2}(\alpha),2}{2}\right)=\left(\frac{N_{E_v/\Bbb
Q_2}(\alpha),2d}{2}\right)=1.$$
\end{proof}

\begin{lem} \label{computation2}
Suppose that $p_i\equiv 1 \mod 8, 1\leq i\leq m$ are distinct
primes. Let $d=p_1p_2\cdots p_m$ such that $2d=r^2+s^2$ with
$r,s\equiv \pm 3 \mod 8$.  If $x_2$ and $y_2$ in $\Q_2$ satisfy
$x_2^2-2dy_2^2=-1$, then the quadratic Hilbert symbol
$$\left(\frac{x_2+y_2\sqrt{2d}, x_0-y_0\sqrt{2d}}{v}\right)=-1$$
where $v$ is the unique place of $E$ above $2$ and $(x_0,y_0)$ is
given as above.
\end{lem}
\begin{proof}
By Lemma \ref{computation}, we have
$$\aligned & \left(\frac{x_2+y_2\sqrt{2d}, x_0-y_0\sqrt{2d}}{v}\right)=\left(\frac{(r-\sqrt{2d})s^{-1}, x_0-y_0\sqrt{2d}}{v}\right) \\
       = \ & \left(\frac{r-\sqrt{2d}, x_0-y_0\sqrt{2d}}{v}\right)\cdot \left(\frac{s, x_0-y_0\sqrt{2d}}{v}\right)  \\
       = \ & \left(\frac{r-\sqrt{2d}, x_0-y_0\sqrt{2d}}{v}\right)\cdot\left(\frac{s,2z_0^2}{2}\right) = -\left(\frac{r-\sqrt{2d}, x_0-y_0\sqrt{2d}}{v}\right). \endaligned  $$
By the Hilbert reciprocity law, one has
$$\left(\frac{r-\sqrt{2d}, x_0-y_0\sqrt{2d}}{v}\right)
=\prod_{\frak p\neq v}\left(\frac{r-\sqrt{2d},
x_0-y_0\sqrt{2d}}{\frak p}\right).$$ Since
$$(r-\sqrt{2d})(r+\sqrt{2d})=-s^2 \ \ \ \text{and} \ \ \
(x_0-y_0\sqrt{2d})(x_0+y_0\sqrt{2d})=2z_0^2$$ with $(r,s)=1$ and
$(x_0,y_0)=1$ respectively, one has $$ord_\frak p(r-\sqrt{2d})\equiv
ord_\frak p(x_0-y_0\sqrt{2d})\equiv 0 \mod 2
$$ for $\frak p\neq v$ and $\frak p<\infty_E$. Since
$x_0-y_0\sqrt{2d}>0$ over $\frak p\in \infty_E$, one obtains that
$$\left(\frac{r-\sqrt{2d}, x_0-y_0\sqrt{2d}}{\frak p}\right)=1 $$ for all $\frak p\neq
v$. One concludes that $$\left(\frac{x_2+y_2\sqrt{2d},
x_0-y_0\sqrt{2d}}{v}\right)=-1.$$
\end{proof}

With a similar argument as in the proof of Theorem \ref{0.1}, we can
prove the following theorem by Lemma \ref{computation} and
\ref{computation2}.
\begin{thm}\label{0.2} Suppose that $p_i\equiv 1 \mod 8, 1\leq i\leq m$ are distinct
primes. Let $d=p_1p_2\cdots p_m$ such that $2d=r^2+s^2$ with
$r,s\equiv \pm 3 \mod 8$. Then the diophantine equation
$x^2-2dy^2=n$ is solvable over $\Z$ if and only if $n$ satisfies the
Artin condition of $H$ and $\Theta$, where $H$ is the Hilbert class
field of $E$ and $\Theta/E$ is a quadratic extension defined as
above.
\end{thm}

\section{Some applications of Theorem \ref{0.2}}
In this section we consider the solvability of the equations
$x^2-2dy^2=-1,\pm 2$ by using Theorem \ref{0.2}. It's well-known
that at most one of the three equations $$x^2-2dy^2=-1,x^2-2dy^2=\pm
2$$ is solvable over $\Z$ (\cite{Per}, p. 106-109).
\begin{lem} \label{computation4} Suppose that $p_i\equiv \pm 1 \mod 8, 1\leq i\leq m$ are distinct
primes. Let $d=p_1p_2\cdots p_m$. If $x_2$ and $y_2$ in $\Q_2$
satisfy $x_2^2-2dy_2^2=2$, then the quadratic Hilbert symbol
$$\left(\frac{x_2+y_2\sqrt{2d}, x_0-y_0\sqrt{2d}}{v}\right)=\begin{cases}1 &\text{ if  }d\equiv 1 \mod
16\\
-1 &\text{ if  }d\equiv 9 \mod 16, \end{cases}$$ where $v$ is the
unique prime of $E$ above $2$  and $(x_0,y_0,z_0)$ is given as in \S
2, $i.e.$ the integers $x_0,y_0$ and $z_0$ are relatively prime with
$x_0>0$ and satisfy $x_0^2-2dy_0^2=2z_0^2$.
\end{lem}
\begin{proof}  The equation $x^2-2y^2=2d$ is solvable over $\Z$.
Choose one solution $(a,b)$ of the equation and obviously $a$ is
even and $b$ is odd. Let $a=2a'$.

First we assume $d\equiv 1 \mod 16$. We will show $b\equiv \pm 1
\mod 8$. Otherwise we have $b\equiv \pm 3 \mod 8$, then we deduce
$a'^2=(d-b^2)/2 \equiv -5 \mod 8$. It is contrary to that $a'\in
\Z$. Similarly we can prove $b\equiv \pm 3 \mod 8$ if $d\equiv 9
\mod 16$. Then we have
$$b\equiv \begin{cases}\pm 1 \mod 8 &\text{ if } d\equiv 1 \mod 16\\ \pm 3 \mod 8 &\text{ if } d\equiv 9 \mod 16\end{cases}.$$

Let $v$ be the unique prime of $E$ above $2$. And let $\xi\in E_v^*
$ with $N_{E_v/\Bbb Q_2}(\xi)=2$. By Lemma \ref{computation}, we
have
$$\aligned & \left(\frac{\xi, x_0-y_0\sqrt{2d}}{v}\right)=\left(\frac{(a-\sqrt{2d})b^{-1}, x_0-y_0\sqrt{2d}}{v}\right) \\
       = & \left(\frac{a-\sqrt{2d}, x_0-y_0\sqrt{2d}}{v}\right)\cdot \left(\frac{b, x_0-y_0\sqrt{2d}}{v}\right)  \\
       = & \left(\frac{a-\sqrt{2d},
       x_0-y_0\sqrt{2d}}{v}\right)\cdot\left(\frac{b,2z_0^2}{2}\right)\\
       =&\begin{cases} \left(\frac{a-\sqrt{2d}, x_0-y_0\sqrt{2d}}{v}\right) &\text{ if } d\equiv 1 \mod 16\\ -\left(\frac{a-\sqrt{2d}, x_0-y_0\sqrt{2d}}{v}\right) &\text{ if } d\equiv 9 \mod 16\end{cases}. \endaligned  $$
Then we only need to show $\left(\frac{a-\sqrt{2d},
x_0-y_0\sqrt{2d}}{v}\right)=1.$ By the Hilbert reciprocity law, one
has
$$\left(\frac{a-\sqrt{2d}, x_0-y_0\sqrt{2d}}{v}\right)
=\prod_{\frak p\neq v}\left(\frac{a-\sqrt{2d},
x_0-y_0\sqrt{2d}}{\frak p}\right).$$ Since
$$(a-\sqrt{2d})(a+\sqrt{2d})=-2b^2 \ \ \ \text{and} \ \ \
(x_0-y_0\sqrt{2d})(x_0+y_0\sqrt{2d})=2z_0^2$$ with $(a,b)=1$ and
$(x_0,y_0,z_0)=1$ respectively, one has $$ord_\frak
p(a-\sqrt{2d})\equiv ord_\frak p(x_0-y_0\sqrt{2d})\equiv 0 \mod 2
$$ for $\frak p\neq v$ and $\frak p<\infty_E$. Since
$x_0-y_0\sqrt{2d}>0$ over $\frak p\in \infty_E$, one obtains that
$$\left(\frac{a-\sqrt{2d}, x_0-y_0\sqrt{2d}}{\frak p}\right)=1 $$ for all $\frak p\neq
v$. One concludes that $$\left(\frac{a-\sqrt{2d},
x_0-y_0\sqrt{2d}}{v}\right)=1.$$
\end{proof}

\begin{prop}\label{9-case}Let $d$ be a positive integer and $d\equiv 9 \mod 16$, then
$x^2-2dy^2=2$ is not solvable over $\Z$.
\end{prop}
\begin{proof} Let $d=p_1^{e_1}\cdots p_m^{e_m}$ where $e_i,1\leq i\leq
m$ are positive integers and $p_i,1\leq i\leq m$ are district
primes. If there is some $p_i$ satisfying $p_i\equiv \pm 3 \mod 8$,
the equation $x^2-2dy^2=2$ is not solvable over $\Z$ since
$\left(\frac{2}{p_i}\right)=-1$. Therefore we can assume $p_i\equiv
\pm 1 \mod 8$ for all $i$.

First we suppose $d$ is square-free. By Lemma \ref{computation4},
one has
$$ \psi_{\Theta/E}( f_E[\prod_{p\leq
\infty}(x_p,y_p)])=\left(\frac{x_2+y_2\sqrt{2d},
x_0-y_0\sqrt{2d}}{v}\right)=-1
$$ for any $\prod_{p\leq \infty}(x_p,y_p)\in \prod_{p\leq \infty}
\bold X_{(2)}(\Bbb Z_p)$, since $\Theta=E(\sqrt{x_0-y_0\sqrt{2d}})$
is unramifed over $E$ everywhere except $v$. Then $x^2-2dy^2=2$ is
not solvable over $\Z$ by Theorem \ref{0.2}.

For general $d=p_1^{e_1}\cdots p_m^{e_m}$ with $p_i\equiv \pm 1 \mod
8$ and $e_i\geq 1$ for all $i$, we can write $d=d'\cdot m^2$ with
$d'$ is square-free and $m\equiv \pm 1 \mod 8$. Then we can see
$d'\equiv d \equiv 9 \mod 16$. Assume $x^2-dy^2=2$ is solvable over
$\Z$. Then $x^2-d'y^2=2$ is solvable over $\Z$, which is contrary to
the above arguments. Therefore we have that $x^2-dy^2=2$ is not
solvable over $\Z$.
\end{proof}

\begin{rem*} If $d$ is a prime and $d\equiv 9\mod 16$ , the unsolvability of $x^2-2dy^2=2$ was proved by Pall (see \cite{Pa}).
\end{rem*}

\begin{lem} \label{computation3} Suppose that $p_i\equiv 1 \mod 8, 1\leq i\leq m$ are distinct
primes. Let $d=p_1p_2\cdots p_m$.  If $x_2$ and $y_2$ in $\Q_2$
satisfy $x_2^2-2dy_2^2=-2$, then the quadratic Hilbert symbol
$$\left(\frac{x_2+y_2\sqrt{2d}, x_0-y_0\sqrt{2d}}{v}\right)=\left(\frac{2}{d}\right)_4$$
where $v$ is the unique prime of $E$ above $2$  and $(x_0,y_0)$ is
given as in \S 2.
\end{lem}
\begin{proof} Since $d=p_1\cdots p_m$ and $p_i\equiv 1 \mod 8,1\leq i \leq m$, we have $x^2+2y^2=2d$ is solvable over $\Z$. Choose
one solution $(a,b)$ of the equation and let $\eta=(a-\sqrt{2d})/b$.
Then $N_{E_v/\Q_2}(\eta)=-2$. By Lemma \ref{computation} we have
$$\left(\frac{x_2+y_2\sqrt{2d}, x_0-y_0\sqrt{2d}}{v}\right)=\left(\frac{\eta,x_0-y_0\sqrt{2d}}{v}\right).$$
So we only need to show
$$\left(\frac{\eta,x_0-y_0\sqrt{2d}}{v}\right)=\left(\frac{2}{d}\right)_4.$$
The Hilbert symbol
$$\aligned \left(\frac{\eta,x_0-y_0\sqrt{2d}}{v}\right)
=\left(\frac{a-\sqrt{2d},x_0-y_0\sqrt{2d}}{v}\right)\cdot
\left(\frac{b,2}{2}\right).
\endaligned  $$
By the Hilbert reciprocity law, we have
$$\left(\frac{a-\sqrt{2d},x_0-y_0\sqrt{2d}}{v}\right)=\prod_{\frak
p\neq v}\left(\frac{a-\sqrt{2d},x_0-y_0\sqrt{2d}}{\frak p}\right).$$
Since
$$(a-\sqrt{2d})(a+\sqrt{2d})=-2b^2 \ \ \ \text{and} \ \ \
(x_0-y_0\sqrt{2d})(x_0+y_0\sqrt{2d})=2z_0^2$$ with $(a,b)=1$ and
$(x_0,y_0,z_0)=1$ respectively, one has
$$ord_\frak p(a-\sqrt{2d})\equiv ord_\frak p(x_0-y_0\sqrt{2d})\equiv
0 \mod 2
$$ for $\frak p\neq v$ and $\frak p<\infty_E$. Since
$x_0-y_0\sqrt{2d}>0$ over $\frak p\in \infty_E$, one obtains that
$$\left(\frac{a-\sqrt{2d}, x_0-y_0\sqrt{2d}}{\frak p}\right)=1 $$ for all $\frak p\neq
v$. So $$\left(\frac{\eta,x_0-y_0\sqrt{2d}}{v}\right)
=\left(\frac{b,2}{2}\right).$$ Since $a^2+2b^2=2d$, we have
$$\left(\frac{d}{l}\right)=1 \text{ for any odd prime } l|a$$
and
$$\left(\frac{2d}{l}\right)=1 \text{ for any odd prime } l|b.$$
Therefore
$$\prod_{p|d}\left(\frac{a}{p}\right)=\prod_{p|d}\left(\frac{a,d}{p}\right)=\left(\frac{a,d}{2}\right)\prod_{\text{odd }l|a}\left(\frac{a,d}{l}\right)
=1\cdot 1=1$$ since $d\equiv 1\mod 8$. And
$$\aligned
\left(\frac{b,2}{2}\right)&=\left(\frac{b,2d}{2}\right)=\prod_{l|b}\left(\frac{b,2d}{l}\right)\cdot\prod_{p|d}\left(\frac{b,2d}{p}\right)\\
&=1\cdot
\prod_{p|d}\left(\frac{b}{p}\right)=\prod_{p|d}\left(\frac{b/a}{p}\right).
\endaligned$$
Since $a^2+2b^2=2d$, we have $(a/b)^2\equiv -2 \mod p \text{ for any
} p|d.$ Hence
$$\left(\frac{a/b}{p}\right)=\left(\frac{-2}{p}\right)_4 \text{ for
any } p|d .$$ Since $p\equiv 1 \mod 8$, then
$\left(\frac{-2}{p}\right)_4=\left(\frac{2}{p}\right)_4 .$ So we
have
$$\left(\frac{\eta,x_0-y_0\sqrt{2d}}{v}\right)=\left(\frac{b,2}{2}\right)=\prod_{p|d}\left(\frac{b/a}{p}\right)=\prod_{p|d}\left(\frac{2}{p}\right)_4=\left(\frac{2}{d}\right)_4.$$

\end{proof}

\begin{prop}\label{aa}  Suppose that $p_i\equiv 1 \mod 8, 1\leq i\leq m$ are distinct
primes. Let $d=p_1p_2\cdots p_m$.  Then

(1) If there exist two integers $r,s\equiv \pm 3\mod 8$ such that
$2d=r^2+s^2$, then the equation $x^2-2dy^2=-1$ is not solvable over
$\Z$.

(2) If the equation $x^2-2dy^2=-2$ is solvable over $\Z$, then
$\left(\frac{2}{d}\right)_4=1$.
\end{prop}

\begin{proof}Denote $\bold X_{(n)}$ to be the affine scheme defined by
$x^2-2dy^2=n$. Let $E=\Q(\sqrt{2d})$ and let $v$ be the unique prime
of $E$ above $2$. Let $\Theta$ and $(x_0,y_0)$ be given as in \S 2.
Then $\Theta/E$ is unramified over all primes except $v$.

(1) By Lemma \ref{computation2}, one has $\left(\frac{\xi,
x_0-y_0\sqrt{2d}}{v}\right)=-1$ for any $\xi\in E_v^* $ with
$N_{E_v/\Bbb Q_2}(\xi)=-1$. This implies that
$$ \psi_{\Theta/E}( f_E[\prod_{p\leq
\infty}(x_p,y_p)])=\left(\frac{x_2+y_2\sqrt{2d},
x_0-y_0\sqrt{2d}}{v}\right)=-1
$$ for any $\prod_{p\leq \infty}(x_p,y_p)\in \prod_{p\leq \infty}
\bold X_{(-1)}(\Bbb Z_p)$. Then $x^2-2dy^2=-1$ is not solvable over
$\Z$ by Theorem \ref{0.1}.

(2) With similar argument as above, the result follows from Lemma
\ref{computation3}.

\end{proof}

If $p\not \equiv 1\mod 8$, the solvability of the three equations
$$x^2-2py^2=-1,x^2-2py^2=\pm 2$$
is well-known (see \cite{D} or \cite{Yo}). If $p\equiv 1 \mod 8$,
the solvability problem is more complicated.
\begin{cor} \label{pall}Let $p$ be an odd prime.

$(1)$ Let $p\equiv 9\mod 16$. If $\left(\frac{2}{p}\right)_4=-1$,
then $x^2-2py^2=-1$ is solvable over $\Z$. If
$\left(\frac{2}{p}\right)_4=1$, then $x^2-2py^2=-2$ is solvable over
$\Z$.

$(2)$ Let $p\equiv 1 \mod 16$. If $\left(\frac{2}{p}\right)_4=-1$,
then $x^2-2py^2=2$ is solvable over $\Z$.
\end{cor}
\begin{proof} Since $p$ is an odd prime, one and only one of the three equations
$$x^2-2py^2=-1,x^2-2py^2=\pm 2$$ is solvable over $\Z$ (\cite{D},
pp 225).

$(1)$ Since $p\equiv 9 \mod 16$, one has the equation $x^2-2py^2=2$
is not solvable over $\Z$ by Proposition \ref{9-case}.

If $\left(\frac{2}{p}\right)_4=-1$, then $x^2-2py^2=-2$ is not
solvable over $\Z$ by proposition \ref{aa}. Therefore $x^2-2py^2=-1$
is solvable.

Suppose $\left(\frac{2}{p}\right)_4=1$. Let $v$ be the unique place
of $E$. Let $\xi_1,\xi_2\in E_v^*$ with $N_{E_v/\Q_2}(\xi_1)=-2$ and
$N_{E_v/\Q_2}(\xi_2)=2$. By Lemma \ref{computation4} and
\ref{computation3}, we have $$\left(\frac{\xi_1,
x_0-y_0\sqrt{2d}}{v}\right)=\left(\frac{2}{p}\right)_4=1 \text{ and
}\left(\frac{\xi_2, x_0-y_0\sqrt{2d}}{v}\right)=-1.$$ By Lemma
\ref{computation}, we have
$$\left(\frac{\xi,
x_0-y_0\sqrt{2d}}{v}\right)=\left(\frac{\xi_1/\xi_2,
x_0-y_0\sqrt{2d}}{v}\right)=1\cdot(-1)=-1$$ for any $\xi\in E_v^*$
with $N_{E_v/\Q_2}(\xi)=-1$, since $N_{E_v/\Q_2}(\xi_1/\xi_2)=-1$.
With the similar arguments as in the proof of Proposition \ref{aa},
one has $x^2-2py^2=-1$ is not solvable over $\Z$. Therefore
$x^2-2py^2=-2$ is solvable over $\Z$.

$(2)$ Suppose $p\equiv 1 \mod 16$. Since
$\left(\frac{2}{p}\right)_4=-1$, one has that $x^2-2py^2=-2$ is not
solvable over $\Z$ by proposition \ref{aa}. Let $v$ be the unique
place of $E$. Let $\xi_1,\xi_2\in E_v^*$ with
$N_{E_v/\Q_2}(\xi_1)=-2$ and $N_{E_v/\Q_2}(\xi_2)=2$. By Lemma
\ref{computation4} and \ref{computation3}, we have
$$\left(\frac{\xi_1,
x_0-y_0\sqrt{2d}}{v}\right)=\left(\frac{2}{p}\right)_4=-1 \text{ and
}\left(\frac{\xi_2, x_0-y_0\sqrt{2d}}{v}\right)=1.$$ By Lemma
\ref{computation}, we have
$$\left(\frac{\xi,
x_0-y_0\sqrt{2d}}{v}\right)=\left(\frac{\xi_1/\xi_2,
x_0-y_0\sqrt{2d}}{v}\right)=(-1)\cdot 1=-1$$ for any $\xi\in E_v^*$
with $N_{E_v/\Q_2}(\xi)=-1$. Therefore $x^2-2py^2=-1$ is not
solvable over $\Z$. Then $x^2-2py^2=2$ is solvable over $\Z$.
\end{proof}
\begin{rem*} The corollary recovers Theorem 3 and 4 in
\cite{Pa}.
\end{rem*}
\bf{Acknowledgment} \it{The work is supported by the Morningside
Center of Mathematics and  NSFC, grant \# 10901150.}


\end{document}